\documentclass[12pt]{amsart}

\usepackage{multicol,fancybox,xspace,fancyhdr,amsmath,amssymb,latexsym,graphicx,amscd}
\newtheorem{proposition}{{\bf Proposition}}[section]
\newtheorem{definition}[proposition]{{\bf Definition}}

\newtheorem{theorem}[proposition]{{\bf Theorem}}
\newtheorem{corollary}[proposition]{{\bf Corollary}}
\newtheorem{remark}[proposition]{{\bf Remark}}
\newtheorem{conj}[proposition]{{\bf Conjecture}}

\usepackage[all]{xy}
\usepackage{epsf}
\usepackage {amssymb,latexsym,amsthm,amscd,amsmath}

\usepackage{multicol,fancybox,xspace,fancyhdr,amsmath,amssymb,latexsym}

\usepackage{graphicx}

\makeindex

\begin{document}

\title{On 3-super bridge knots}
\author{Rama Mishra}

\address{Indian Institute of Science Education and
Research, Pune 411008, India. \newline \mbox{} \hspace{1.5mm}
{\small \em r.mishra@iiserpune.ac.in}}

\date{}

%%%%%%%%%%%%%%%%%%%%%%%%%%%%%%%%%%%%%%%%%%%%%%%%%%%%%%%
      % This is the title page.

%%%%%%%%%%%%%%%%%%%%%%%%%%%%%%%%%%%%%%%%%%%%%%%%%%%%%%%
%\include{declare}    % This is to be attached separately.
%\pagebreak
%%%%%%%%%%%%%%%%%%%%%%%%%%%%%%%%%%%%%%%%%%%%%%%%%%%%%%%%
%%%%%%%%%%%%%%%%%%%%%%%%%%%%%%%%%%%%

%\thispagestyle{empty}

\begin{abstract}

It is known that there are only finitely many knots with
super bridge index 3. Jin and Jeon  have provided a list of
possible such candidates. However, they conjectured that the only
knots with super bridge index 3 are trefoil and the figure eight
knot. In this paper, we prove that the $5_2$ knot and the $6_2$ knot are also 3-super
bridge knots by providing a polynomial representation of these knots in degree $6.$
This also answers a question asked by Durfee and O'Shea in
their paper on polynomial knots: is there any 5-crossing knot in
degree 6?

\end{abstract}

\maketitle



%\end{frontmatter}

\section{Introduction}

In the recent years mathematical theory of knots has become an
interdisciplinary subject. In the physical situations we come
across specific conformation of knots and we have to draw
inferences based on that. Knot invariants are the main tools to
use knot theory anywhere. Thus, to apply knot theory in physical
scenario we need to use the knot invariants which are dependent on
the conformation itself. Some of such important knot invariants
are the {\it crossing number,  bridge number} and the {\it
unknotting number}. The idea behind defining each one of them is
similar and is based on the following theme: Fix a conformation of
a knot in the space. Define this quantity as minimum  for this
conformation and then minimize it over all conformation of knots.
In \cite{kui} Nicolas Kuiper defined a new invariant called {\it
Super bridge index} for knots. It differs from the bridge index in
the following sense: for a fixed conformation we consider the
maximum (instead of minimum as in case of bridge index) number of
local maxima coming from the projection to an axis. We then
minimize it over all conformations. This is a physical knot
invariant as it depends on how conformations of the knots sit in
Euclidean 3-space. For a given knot $K$ the bridge index $b[K]$
and the super bridge index $sb[K]$ are related as: $b[K]\leq
sb[K].$ Later Kuiper proved that $b[K]<sb[K]$ for any non trivial
knot $K.$. Thus for a non trivial knot $K$, the super bridge index
$sb[K]$ has to be at least $3$. Kuiper  computed the super bridge
index for all torus knots of type $(p,q)$ and prove that it is $min\{2p, q\}$.
As a direct corollary to the Kuiper's result on torus knots one obtains that there are
infinitely many knots with super bridge index $n$ for each even integer $n\geq
4.$ It is believed that there are infinitely many $n-$ super bridge knots for every integer $n\geq 4.$
The situation regarding knots with super bridge index $3$
turned out more interesting. They are certainly 2 bridge knots.
Using the concept of quadrisecants, Jin and Jeon \cite{jj1} proved
that there are only finitely many knots with super bridge index 3.
They in fact, made a list of possible knots with super bridge
index 3 to be: $3_1, 4_1, 5_2, 6_1, 6_2, 6_3,
 7_2, 7_3, 7_4, 8_4$ and $8_9.$  Among these, only for $3_1$ and $4_1$
 it has been proved that the super bridge index is $3$. It was conjectured that
 these two knots are perhaps the only knots with super bridge index 3.

As the super bridge index is dependent on conformation of knots,
various upper bounds have been found using some of the knot
parameterizations in terms of stick number \cite{rand}, Harmonic degree
\cite{trau} and polynomial degree \cite{ad}.  We have used the approach
of polynomial parametrization. A positive integer $d$ is said to
be the polynomial degree of a knot $K$ if one can find an
embedding $t\to (f(t),g(t),h(t))$ from $\mathbb{R}$ to
$\mathbb{R}^3$ where $f(t),g(t)$ and $h(t)$ are polynomials over
$\mathbb{R}$ with $deg(f(t)<deg(g(t))<deg(h(t))=d$ such that the
one point compactification of the image of this embedding is
equivalent to the knot $K$ and $d$ is the least positive integer
with this property. If a knot $K$ has polynomial degree $d$ then
it can be proved that its super bridge index $sb[K]\leq
\frac{d+1}{2}.$ Given a knot type $K$ determining its polynomial
degree is an interesting question. Trefoil knot and the figure
eight knot are known to have polynomial degree $5$ and $6$
respectively. It is easy to prove that trefoil is the only knot
that can be represented in degree $5$. Durfee and O'Shea had
raised a question: are there any five crossing knot with
polynomial degree $6$? Using Kuiper's result on super bridge index
of torus knots the knot $5_1$ cannot have polynomial degree $6$.
However the $5_2$ knot was still the hope. In this paper we
produced a polynomial parametrization of $5_2$ knot in degree $6$.
This also proved that the super bridge index of $5_2$ knot must be
$3$ and hence disproves the conjecture of Jin and Jeon. We also found a degree $6$ polynomial representation for $6_2$ knot which proves that its super bridge index must be $3.$

This paper is organized as follows: In section 2, we define the
super bridge index of a knot and provide an overview of the
results that are relevant for the paper. We discuss about
polynomial knots and its connection with the super bridge index in
section 3. In section 4, we provide a degree 6 polynomial
parameterizations of the $5_2$ and the $6_2$ knots and hence prove that their super
bridge index must be 3.

\section{The Super bridge index of a knot}

 An ambient isotopy class of  smooth embedding of the unit circle
 $S^1$ in three space $\mathbb{R}^3$ is called a {\em knot type} or
 simply a knot. An embedding $\varphi: S^1\to \mathbb{R}^3$ is always identified with
 the image $\varphi(S^1)\subset \mathbb{R}^3.$ Thus we denote a knot by $[K]$. In an isotopy
 class, image of each embedding may look different but they are
 equivalent in the sense that one can be transformed into the other
 using smooth deformation. An individual embedding $\varphi$ in the isotopy
 class $[K]$ is called a {\em knot conformation}. Thus one can say that
 a knot is an equivalence class of knot conformations.

 Let $\varphi: S^1\to \mathbb{R}^3\in [K]$ be a knot conformation of $[K]$. Let
 $\overrightarrow{v}$ be a vector in $\mathbb{R}^3$ such that
 $\|\overrightarrow{v}\|=1.$ Let $\Pi_v:\mathbb{R}^3\to
 \mathbb{R}^2$ be the projection map onto the plane perpendicular
 to $\overrightarrow{v}.$ The number of local maxima in the
 projected curve $\Pi_v\circ \varphi$ is finite. Let $m_{\varphi}$ and
 $M_{\varphi}$ denote the minimum and the maximum number of local
 maxima in $\Pi_v\circ \varphi$ taken over all $\overrightarrow{v}\in \mathbb{R}^3$ such that
 $\|\overrightarrow{v}\|=1.$

In \cite{schu} Schubert defined the {\em bridge index $b([K])$} of a
knot $[K]$ to be the $min_{\varphi\in[K]}\{m_{\varphi}\}$. He
proved that for any non trivial knot $[K]$ its bridge index
$b([K])$ has to be at least $2.$ A knot with bridge index $n$ is
referred as an $n-$bridge knots. Schubert completely classified
all $2-$ bridge knots and also determined the bridge index of all
torus knots. Many mathematicians have worked on the problem of
computing the bridge index and right now the bridge index is known
for all knots up to 10 crossings.

In \cite{kui} Nicolas Kuiper came up with a new idea. Since any knot conformation is compact.
If we consider the number of local maxima of the projection into the plane perpendicular to each direction where it is a
regular projection,
there will be a maximum value. He came up with the following\\

\begin{definition} {\rm  Let $\varphi: S^1\to
\mathbb{R}^3\in [K]$ be a knot conformation of $[K]$. Then
$M_{\varphi}$ the maximum number of local maxima in $\Pi_v\circ
\varphi$ taken over all $\overrightarrow{v}\in \mathbb{R}^3$ such
that
 $\|\overrightarrow{v}\|=1$ and the projection is a regular projection is called the {\em Super bridge number} of the
 knot conformation
 $\varphi$ and is denoted by $Sb(\varphi)$.}
\end{definition}

\begin{definition}
%\noindent{\bf\large Definition 2.2.} The {\em super bridge index}
{\rm $Sb([K])$  of a knot $[K]$ is defined as:
$$Sb([K]):=min_{\varphi\in [K]}\{Sb(\varphi)\}.$$}
\end{definition}

From the definition it is clear that for any knot {\nolinebreak $[K]$
$b([K])\leq Sb([K])$.} In his paper \cite{kui} Kuiper proved the
following results:

\begin{enumerate}
\item {\em If $[K]$ is any non trivial knot then $b([K])<
Sb([K]).$ Thus for any non trivial knot $[K]$ its super bridge
index $Sb([K])\geq 3.$}

\item  {\em Let $[K_{p,q}]$ denote a torus knot of type $(p,q)$
with $p<q$. Then $Sb([K_{p,q}])=min\{2p,q\}.$}
\end{enumerate}

A knot $K$ whose super bridge index is $n$ is referred as {\em $n-$
super bridge knot}. From (2), it is clear that for each even
integer $n\geq 4$ there are infinitely many $n-$super bridge
knots.  In fact, it is conjectured \cite{jj1} that not just even but
for all integers $n\geq 4$ there are infinitely many $n-$super
bridge knots. The case for $n=3$ was studied separately. In
\cite{jj1} Jin and Jeon  proved the following theorems.

\begin{theorem} {\em There are only finitely
many 3-super bridge knots}.
\end{theorem}
\begin{theorem} {\em The only possible 3-super
bridge knots are:\\
 $3_1, 4_1, 5_2, 6_1, 6_2,6_3,7_2, 7_3, 7_4,
8_4, 8_7$ and $8_9.$}.
\end{theorem}

Later Adams and his research group \cite{adam}  proved that the knot $8_7$ has super bridge
index $4.$ Thus the list cuts down to only $11$ possible knots. However among the possible candidates in the list
of 3-super bridge knots only $3_1$ and $4_1$ are confirmed to have
super bridge index $3.$  In \cite{jj1} Jin made the following:\\

\begin{conj} $3_1$ and $4_1$ are the only
$3-$super bridge knots.
\end{conj}

In section 4, we will prove that $5_2$ is a 3-super bridge knot
and hence disprove the above conjecture.

\section{Polynomial knots}

\begin{definition} {\rm A non compact knot or a long
knot $K$ is a proper, smooth embedding of $\mathbb{R}$ in
$\mathbb{R}^3$ which is monotone outside a closed interval.}
\end{definition}

It is clear that any non compact knot will have a unique extension
as a smooth embedding of $S^1$ in $S^3$
which is our classical knot.\\

\begin{definition} {\rm Two long knots $k_0$ and
$k_1$ are said to be topologically equivalent if their extensions
as classical knots are ambient isotopic. Or equivalently if there
exists an orientation preserving diffeomorphisms
$h:\mathbb{R}^3\to \mathbb{R}^3$ such that $h(k_0)=k_1.$}
\end{definition}

A typical long knot is an embedding $t\to (\varphi_1(t),
\varphi_2(t), \varphi_3(t))$ where $\varphi_i$ are smooth real
valued functions such that for $|t|\to \infty$ each
$|\varphi_i(t)|\to \infty.$

\begin{definition} {\rm A long knot defined by an
embedding of the form $t\to (f(t), g(t), h(t))$, where $f(t),g(t)$
and $h(t)$ are real polynomials, is called a polynomial knot.}
\end{definition}

It has been proved (\cite{ars} and \cite{rs}) that each long knot
is topologically equivalent to some polynomial knot.

\begin{definition} {\rm A polynomial knot defined by
$t\to (f(t), g(t), h(t))$ is said to have degree $d$ if
$deg(f(t))=deg(g(t))=deg(h(t)) =d.$}
\end{definition}

Since the polynomial automorphisms of the form $(x,y,z)\to (x-\alpha z, y,z)$ and $(x,y,z)\to (x, y-\beta z,z)$,
for $\alpha$ and $\beta$ reals, are orientation preserving, it leads to a polynomial knot to acquire the form
$t\to (f(t),g(t),h(t))$ where $deg(f(t))<deg(g(t))<deg(h(t))$ and none of these degrees lie in the semi group generated
by the other two. Suppose
$t\to (f(t),g(t),h(t))$ is a polynomial knot of this form and $deg(f(t))=d_1,$ $deg(g(t))=d_2$ and $deg(h(t))=d_3,$
the triple $(d_1,d_2,d_3)$ is called a {\it degree sequence} of the polynomial knot. Also if we have a polynomial knot
$K$ given by
$t\to (f(t),g(t),h(t))$ with $deg(f(t))=d_1,$ $deg(g(t))=d_2$ and $deg(h(t))=d_3,$ such that $d_1<d_2<d_3$ and none of
them lie in the semi group generated by the other two. Then we can choose $\epsilon>0$ and $\delta>0$ such that the
polynomial knot given by $t\to (f(t)+\epsilon t^{d_3},g(t)+\delta t^{d_3} ,h(t))$ is topologically equivalent to $K.$\\
\label{}

Thus, we can rephrase the definition 3.4 as:  A polynomial knot
defined by $t\to (f(t), g(t), h(t))$ is said to have degree $d$ if
$deg(f(t))<deg(g(t))<deg(h(t))$ and $deg(h(t)) =d.$

It is easy to note that if a polynomial knot $K$ has degree $d$, we can obtain polynomial knots of degree $d+k$ for
each $k\geq 1$ which are topologically equivalent to $K.$\\

\begin{definition} {\rm A positive integer $d$ is
said to be the minimal degree for a knot $K$ if there is a
polynomial knot defined by $t\to (f(t), g(t), h(t))$ which is
topologically equivalent to $K$ with
$deg(f(t))<deg(g(t))<deg(h(t))$ and $deg(h(t)) =d$ and no polynomial knot with degree less than $d$ is equivalent to $K$.}
\end{definition}

\begin{remark}
\begin{itemize}
\item{(i)} Clearly the minimal degree of a knot is a knot
invariant. \item{(ii)} If a knot $K$ is represented by $t\to
(f(t), g(t), h(t))$ then the polynomial knot given by $t\to (f(t),
g(t), -h(t))$ represents the mirror image of $K$. \item{(iii)}
Thus a knot and its mirror image have same minimal degree. Hence
the minimal degree cannot detect the {\it chirality} in a knot.
\end{itemize}
\end{remark}

Certain numerical knot invariants can be inferred from the degree
of a polynomial knot. In this connection we have:

\smallskip

\begin{theorem}  Let $K$ be a knot with a polynomial representation in degree
$d.$ Let $C(K)$, $b(K)$ and $S(K)$ respectively denote the minimal crossing number, the bridge index and the super bridge index of $K$. Then we have the following estimates
\end{theorem}

\begin{enumerate}

\item $c(K)\leq \frac{(d-2)(d-3)}{2}.$

\item  $ b(K)\leq \frac{(d-1)}{2}.$

\item $ s(K)\leq \frac{(d+1)}{2}.$

\end{enumerate}

\noindent{\bf Proof.}  (1). We can assume that the knot $K$ has a polynomial representation $t\to(f(t),g(t),h(t))$ with
$\deg(f(t))\leq d-2,$ $\deg(g(t))\leq d-1$ and $\deg(h(t))=d$. The number of double points for any plane curve $(x(t),y(t))$ can be obtained by finding the 
roots of the resultant polynomial of $\frac{x(t)-x(s)}{t-s}$ and $\frac{y(t)-y(s)}{t-s}$ with $t\neq s$, in one of the variables $t$ or $s$. The double points 
will be half of the number of the roots of this resultant polynomial. In this case the minimum number of double points is obtained on the $xy$ plane. Here the 
resultant will be a polynomial of degree $(d-3)(d-2)$ and will have at most $(d-3)(d-2)$ roots. Thus $C(K)\leq \frac{(d-2)(d-3)}{2}.$

\smallskip

(2). The number of local maxima in any direction is at most half of the roots of the derivative of the function in that direction. In this case the minimum number of local maxima among all the regular projections is less than or equal to $\frac{d-3}{2}.$ Thus when we take the one point compactification of this polynomial knot, the minimum number of local maxima is at most $\frac{d-3}{2}+ 1= \frac{d-1}{2}.$ Hence $b(K)\leq \frac{d-1}{2}.$

\smallskip

(3). The maximum number of local maxima among all the regular projections is $\frac{d-1}{2}.$ After taking the one point compactification, the maximum number 
of 
local maxima for this knot is at most $\frac{d-1}{2}+1=\frac{d+1}{2}.$ Hence $s(K) \leq \frac{d+1}{2}.$

This completes the proof of the theorem. \hfill$\Box$

\section { The $5_2$ knot}

We show that there exists a polynomial representation of $5_2$ knot in degree $6.$ We prove the following

\begin{theorem} The knot $5_2$ can be represented by a polynomial embedding $t\to (f(t),g(t),h(t))$ where $\deg(f(t))=4,$
$\deg(g(t))=5$ and $\deg(h(t))=6.$
\end{theorem}

\smallskip
\noindent{\bf Proof.}  Consider the plane curve
$$ (x(t),y(t) = ((t - 2) (t + 4) (t^2 - 11),t (t^2 - 6) (t^2 - 16)).$$ This can be verified that $x^\prime(t)$ and $y^\prime(t)$
do not have any common zero. Also, the self intersections of this curve can be determined by finding the common zeros of

$$ X(t,s)=\frac{x(t)-x(s)}{t-s}=0\:\:{\rm and}\:\: Y(t,s)=\frac{y(t)-y(s)}{t-s}=0$$ Here $t\neq s.$

 After performing the computations in {\em mathematicae} We find that the double points (approximately) occur at the pair of parametric
values $(t_1,s_1)= (-4.21,3.43)$,$(t_2,s_2)= (-3.85,2.08)$,
$(t_3,s_3)= (-3.01,1.79)$, $(t_4,s_4)=(-2.05,3.84)$,and
$(t_5,s_5)= (.105, 4.03)$. Clearly this is a regular projection
of a knot on $xy$ plane with at most $5$ crossings. Its graph is given in Figure 1.

\begin{center}
{\includegraphics[height= 2in]{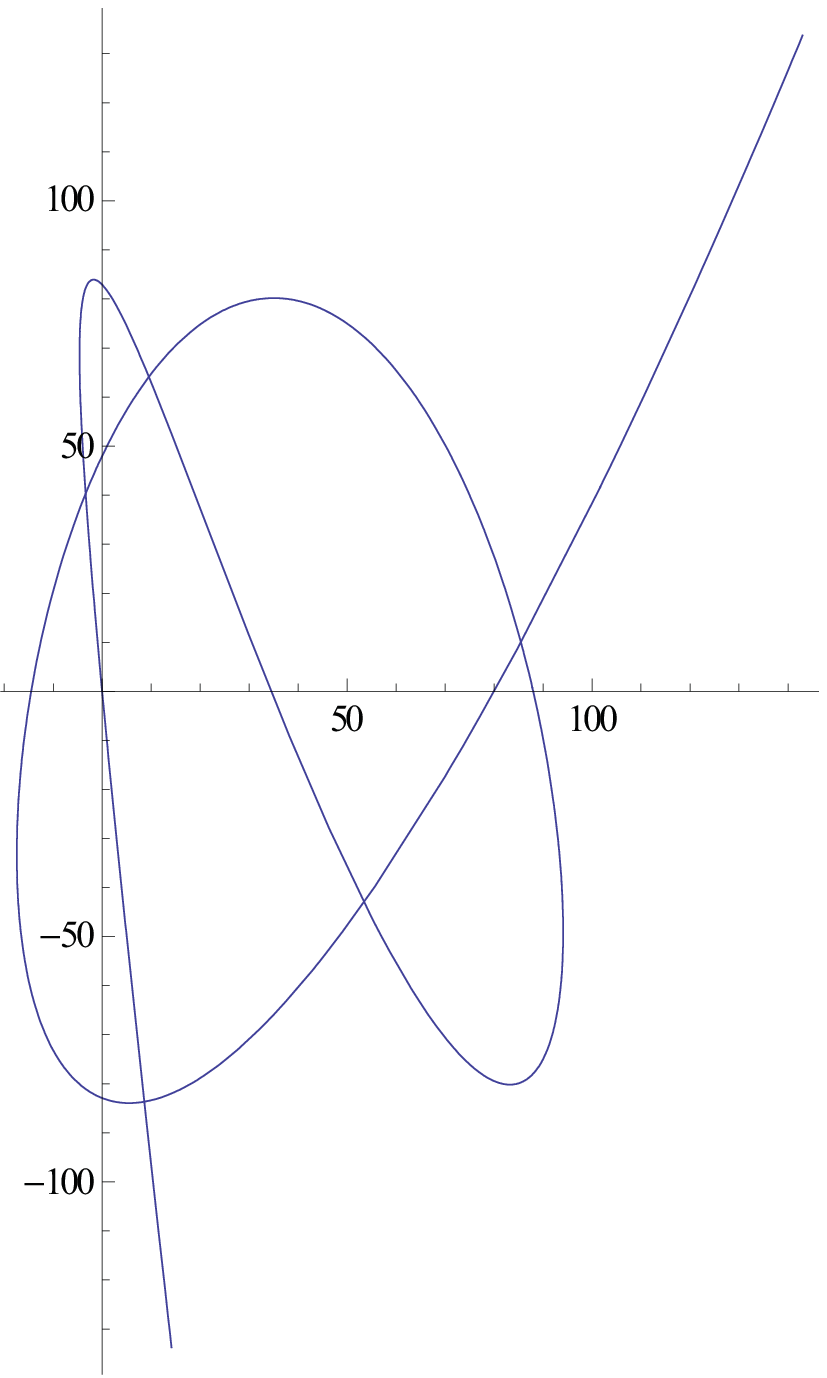}}

  Figure 1

 \end{center}

Now we need to find a polynomial $h(t)$ such that $h(t_i)-h(s_i)$ is negative for $i$ odd and positive for $i$ even.
Let  $$h(t)=t^6+At^5+Bt^4+Ct^3+Dt^2+Et$$ be a degree $6$ polynomial. We want to find the values of the coefficients
$A,$ $B$, $C$, $D$ and $E$ such that $h(t_i)-h(s_i)= -100$ for $i$ odd and $h(t_i)-h(s_i)= 100$ for $i$ even.

 This gives us a system of $5$ linear equations in $5$ variables $A,B,C,D$ and $E.$ This system has a unique solution because the determinant of the coefficient matrix is $5123.92$ which is non-zero.   In this example, {\it mathematica} computations give

 $A=-1505.83,$ $B=-293.032,$ $C=32625.7,$ $D=5323.59$ and $E=-138788.$  Thus we have $z(t)=t^6-1505.83t^5-293.032t^4+32625.7t^3+5323.59t^2-138788t$.
 We can clearly check the values of $h(t)$ at parametric values giving us the double point. $h(t_1)=149216$ and $h(s_1)=149320$.
 Thus $h(t_1)<h(s_1)$, i.e., $(t_1,s_1)$ is an under crossing. Similarly $h(t_2)=-35983.8>-36081.5=h(s_2)$ means $(t_2,s_2)$ is an over crossing,
 $h(t_3)=-75002<-74901=h(s_3)$, means $(t_3,s_3)$ is an under crossing, $h(t_4)=75231>75133=h(s_4)$, i.e., $(t_4,s_4)$ is an over crossing and at
 last $h(t_5)=-14476<-11152=h(s_5)$ proving that $(t_5,s_5)$ is an under crossing.

Thus we have a polynomial parametrization of a knot in degree $6$:

\begin{eqnarray*}
x(t)&=&(t - 2) (t + 4) (t^2 - 12)\\
y(t)&=& t (t^2 - 6) (t^2 - 16)\\
z(t)&=& t^6-1505.83t^5-293.032t^4+32625.7t^3+5323.59t^2-138788t
\end{eqnarray*}
A $3d$ plot of this embedding must look like as shown in Figure 2. 

\begin{center}
{\includegraphics[height=1.5in]{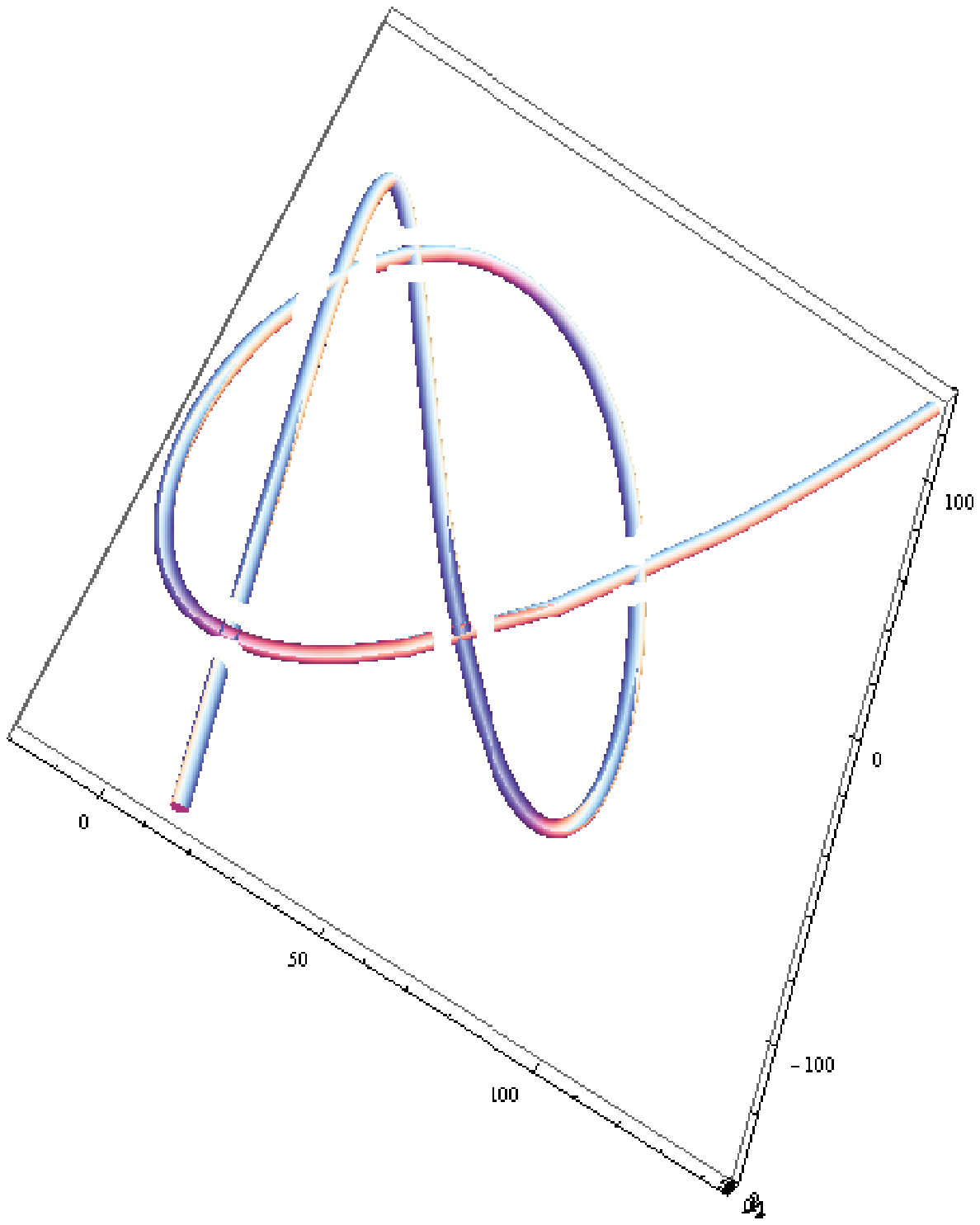}}

  Figure 2

 \end{center}
 
 The one point compactification of this embedding will represent a knot diagram as shown in Figure 3. We would like to identify the knot type of this knot diagram. For this we compute the Jones polynomial of this diagram. Let $D$ denote the diagram shown in Figure 3. The calculations for the Kauffman bracket of the diagram $D$ are shown in Figure 4.

\begin{center}
{\includegraphics[height=2in]{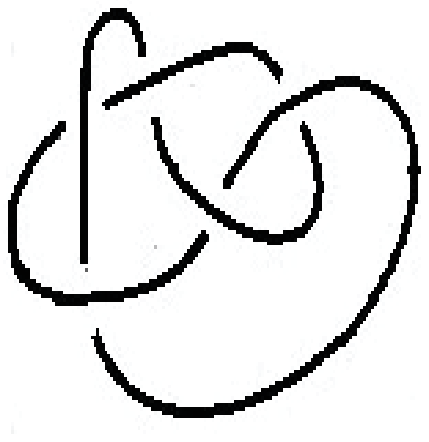}}
 \end{center}
 
\begin{center}
Figure 3
\end{center}

 \begin{center}
 {\includegraphics[height= 3in]{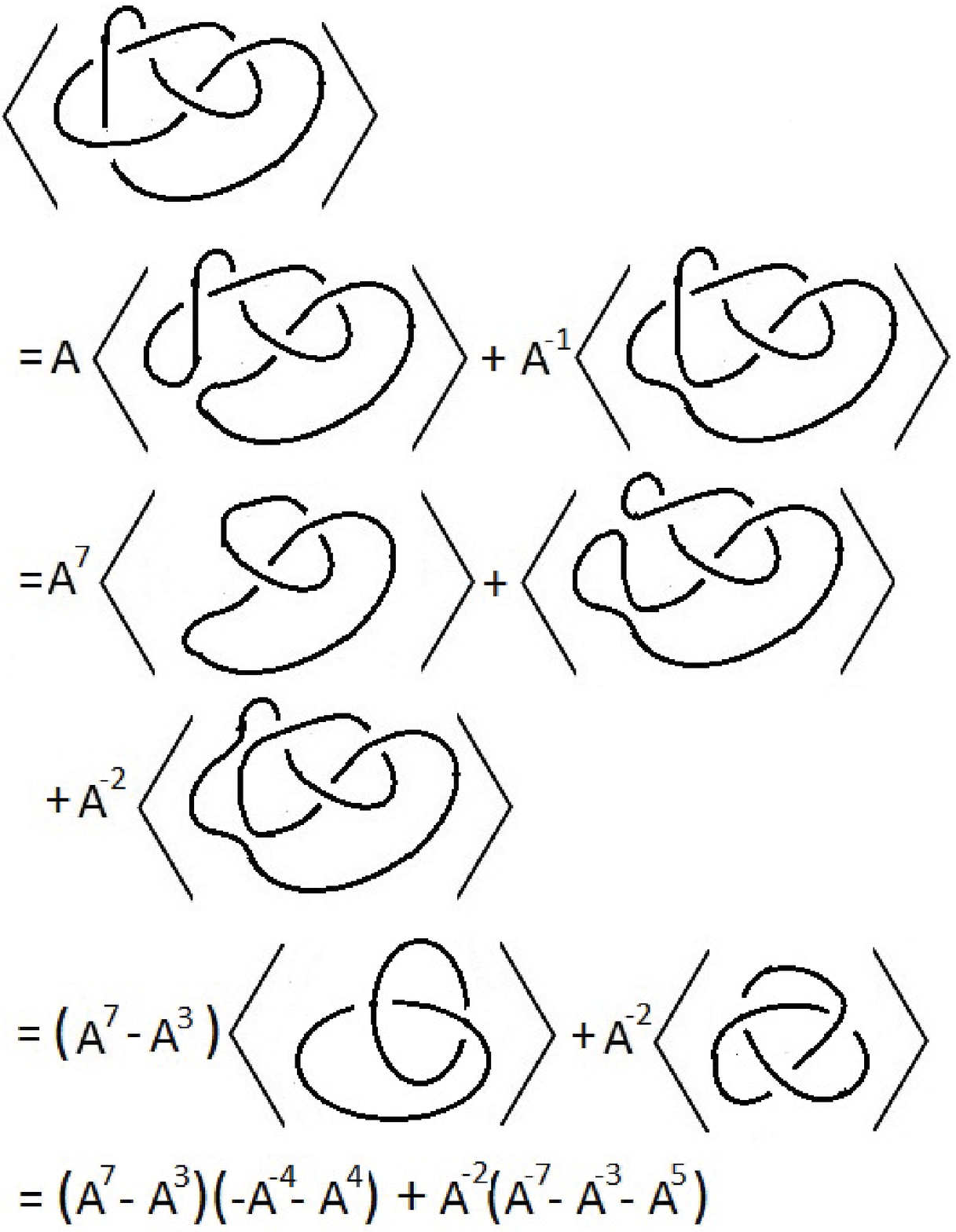}}

  Figure 4
 \end{center}
 
 \medskip
 
 After simplification we obtain
 
\begin{eqnarray*}
\langle D\rangle &=&(A^7-A^3)(-A^{-4}-A^4)+A^{-2}(A^{-7}-A^{-3}-A^5)\\
                 &=& -A^3-A^{11}+A^{-1}-A^7+A^{-9}-A^{-5}-A^3\\
                 &=& -A^{11}+A^7-2A^3+A^{-1}-A^{-5}+A^{-9}\\
   \end{eqnarray*}

Thus the Kauffman polynomial of the knot $K$ which is given by the knot diagram $D$ is,

\begin{eqnarray*}
           f_K(A)&=&(-A)^{-3w(D)}\langle D\rangle \\
                 &=& (-A)^{(-3)(5)}(-A^{11}+A^7-2A^3+A^{-1}-A^{-5}+A^{-9})\\
                 &=& -A^{-15}(-A^{11}+A^7-2A^3+A^{-1}-A^{-5}+A^{-9})\\
                 &=& A^{-4}-A^{-8}+2A^{-12}-A^{-16}+A^{-20}-A^{-24}\\
\end{eqnarray*}

Thus the Jones polynomial $V_K(q)$ is obtained by making the substitution $A^{-4}=q$ and we get
$$V_K(q)= q-q^2+2q^3-q^4+q^5-q^6.$$  This is the Jones polynomial of $5_2$ knot. It is known that up to $9$ crossings all knots have distinct Jones polynomial( \cite{kauf}) and hence it proves that the diagram $D$ indeed represents the $5_2$ knot. Thus we have a degree $6$ polynomial representation of $5_2$ knot.
This completes the proof of the theorem. \hfill$\Box$

 From this theorem we have the following

\begin{corollary} 
The minimal degree of $5_2$ is 6.
\end{corollary}
\smallskip

\noindent{\bf Proof.} Follows from $(1)$ of Theorem 3.7.

\begin{corollary} 
The super bridge index of $5_2$ is 3.
\end{corollary}

\smallskip

\noindent{\bf Proof.} Follows from $(3)$ of Theorem 3.7.

\section{The $6_2$ knot}
 
\begin{theorem} The knot $6_2$ can be represented by a polynomial embedding $t\to (f(t),g(t),h(t))$ where $\deg(f(t))=4,$
$\deg(g(t))=5$ and $\deg(h(t))=6.$
\end{theorem}

\smallskip

\noindent{\bf Proof.}  Consider the plane curve
$$ (x(t),y(t) = ( t^4 - 27 t^2 + t, t^5 - 36 t^3 + 260 t).$$ This can be verified that $x^\prime(t)$ and $y^\prime(t)$
do not have any common zero. After performing the computations in {\em mathematicae} we find that the double points  occur (approximately) at the pair of parametric
values $(t_1,s_1)= (-5.201,-.363),\:\:\:\:(t_2,s_2)= (-5.118,5.078)\:\:\:\:,
 (t_3,s_3)= (-4.698,2.31),\:\:\:\: (t_4,s_4)=(-3.090,3.233),\:\:\:\:
(t_5,s_5)= (-2.226, 4.651)$ and $(t_6,s_6)=(.252, 5.172).$ Clearly this is a regular projection
of a knot on $xy$ plane with at most $6$ crossings. Its graph is given in Figure 5. 
\begin{center}
{\includegraphics[height= 2.5in]{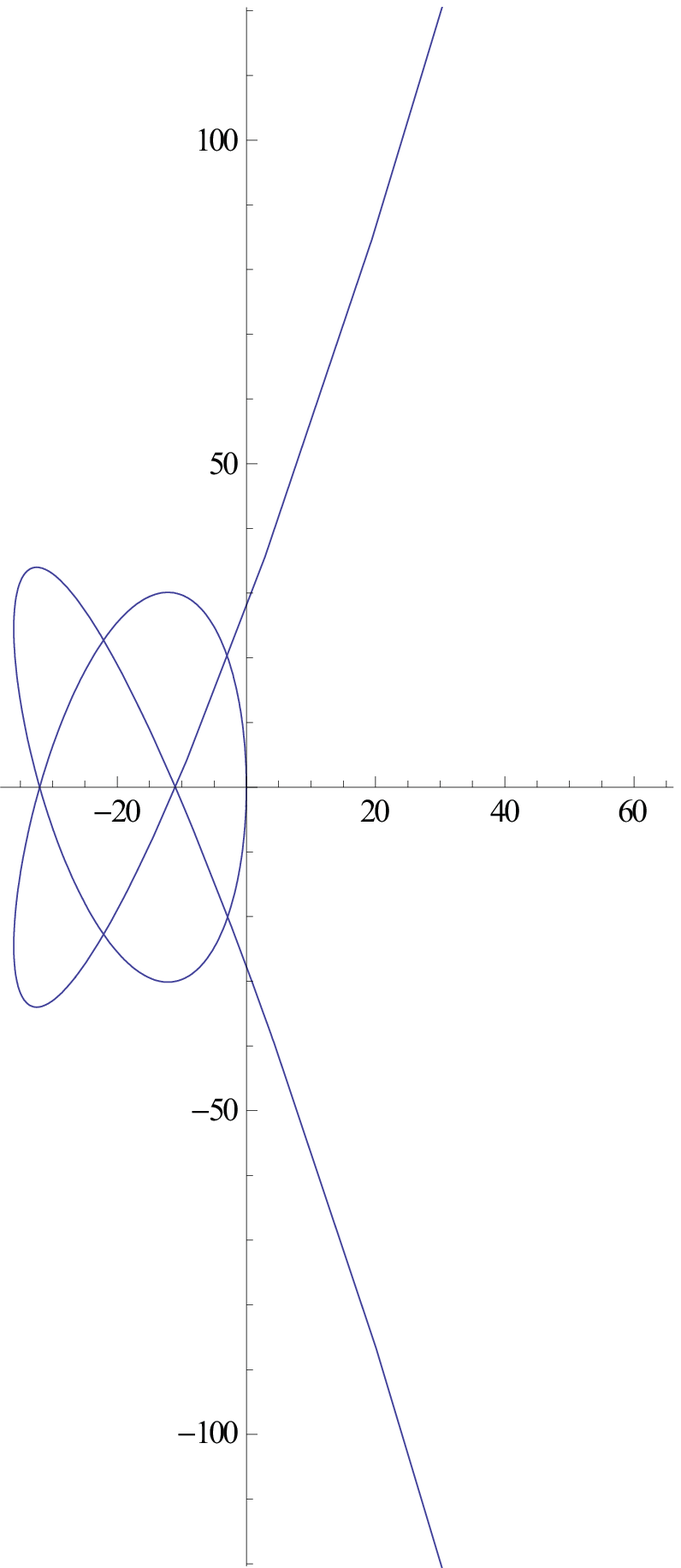}}

  Figure 5

 \end{center}
 
 To provide an over/under crossing information, let $h(t)=At^6+Bt^5+Ct^4+Dt^3+Et^2+Ft$ be a degree $6$ polynomial such that $h(t_i)-h(s_i) = -100$ for $i$ odd and $h(t_i)-h(s_i) = 100$ for $i$ even. This will give us a system of $6$ linear equations in $6$ variables. The determinant of the coefficient matrix of this system is found (using mathematica) to be $-5.22794\times 10^6$ which is non zero. The values of the unknowns are estimated as :
$A=-0.0221563$, $B=-413.2$, $C=3202.02$, $D=14878.7$, $E=-86446.7$ and $F=-104260.$ Thus
$$h(t)= -.02215t^6 - 413.2t^5 + 3202.02t^4 + 14878.7t^3 - 86446.7t^2 - 
 104260 t.$$

 At each double point we check that $h(t_1)=25638.2,\:\:h(s_1)=25801.5,\:\: h(t_2)=-76678.4,\:\: h(s_2)= -76756.5,
h(t_3)=-455081,\:\: h(s_3)=-454735,\:\:\:\: h(t_4)=-533917,\:\:\:\:\:\:$ $h(s_4)=\:\:-534002,\:\: h(t_5)=-259180,\:\: h(s_5)=-259134, h(t_6)= -31512.6,\:\: h(s_6)=-31604.6.$ Thus $(t_i,s_i)$  are under crossing points for $i=1,3,5$ and are over crossing points for $i=2,4,6.$ Thus the knot diagram of the embedding $t\to (( t^4 - 27 t^2 + t, t^5 - 36 t^3 + 260 t,-.02215t^6 - 413.2t^5 + 3202.02t^4 + 14878.7t^3 - 86446.7t^2 -  104260 t )$ will look as is shown in Figure 6.

\begin{center}
{\includegraphics[height=2in]{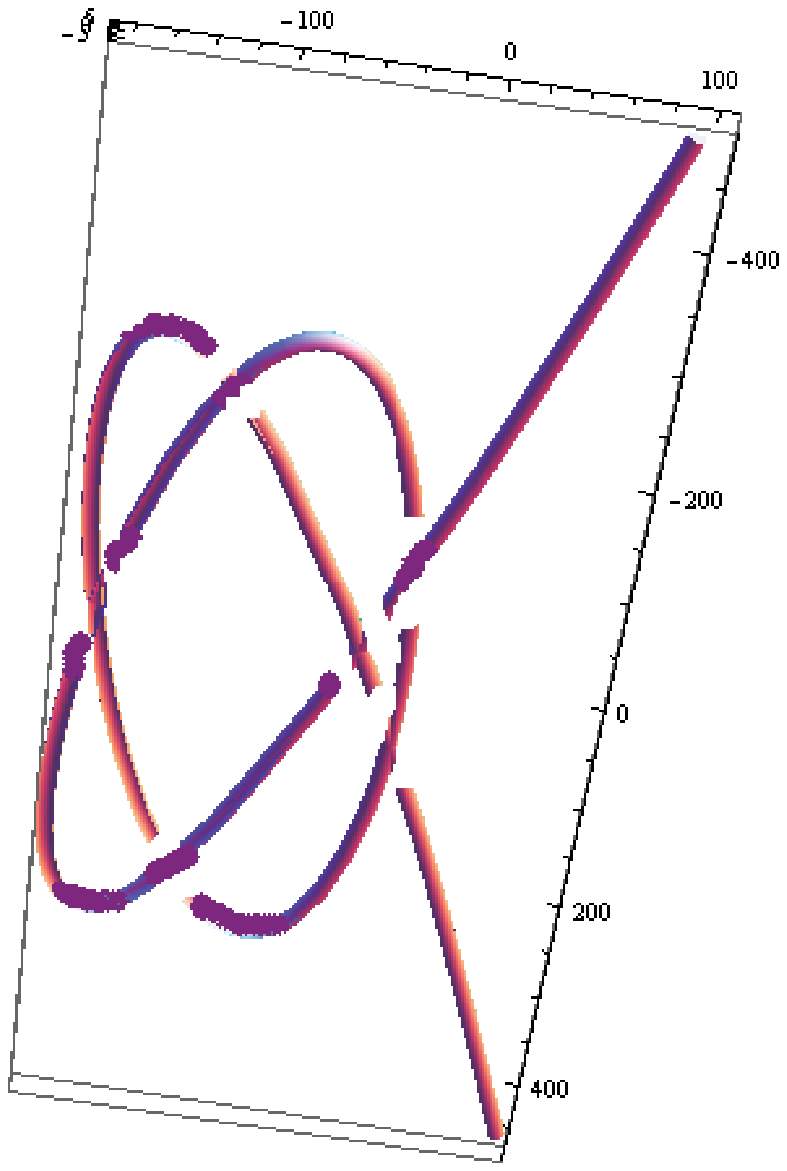}}

  Figure 6

 \end{center}

After taking the one point compactification, this diagram will appear as in Figure 7.

\begin{center}
{\includegraphics[height= 2in]{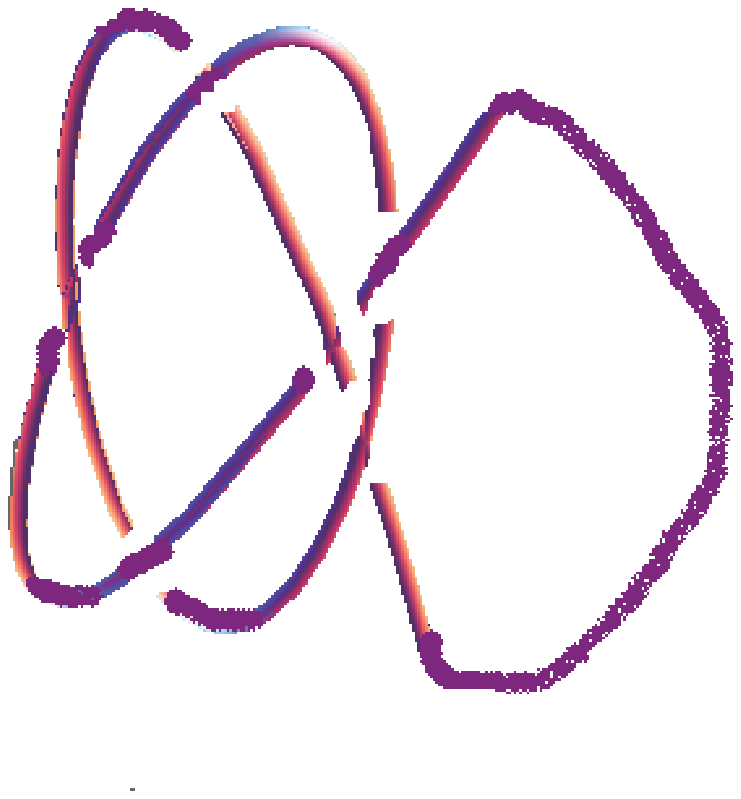}}

  Figure 7

 \end{center}

This is a knot diagram of a knot with at most $6$ crossings. We compute the Jones polynomial of this knot diagram and find that it is $ q-1+2q^{-1}-2q^{-2}+2q^{-3}-2q^{-4}+q^{-5}.$ This is the Jones polynomial of the $6_2$ knot. Using a similar argument as in the previous example we conclude that the diagram shown in Figure 7 represents the $6_2$ knot.
Thus we have found a degree $6$ polynomial representation of the $6_2$ knot. This completes the proof. \hfill$\Box$.

\medskip

\begin{corollary}
The super bridge index of $6_2$ knot is $3.$
\end{corollary}

\noindent{\bf Proof.}  Follows from $(3)$ Theorem 3.7.


\begin{thebibliography}{99}
\bibitem{adam}
Colin Adams and 2007 SMALL Research Group, Superbridge Number of Knots,
preprint, 2007.

\bibitem{ad} Alan Durfee and Don Oshea, Polynomial knots, http://arxiv.org/pdf/math/0612803v1.pdf.
\bibitem{jj1}  C.B. Jeon and G.T. Jin, There are only finitely many 3-superbridge knots, J. Knot Theory Ramifications Vol.10, 2, 2001, 331-343.
\bibitem{jj2} C.B. Jeon and G.T. Jin, A computation of superbridge index of knots, J. Knot Theory Ramifications Vol.11, 3, 2002, 461-473.
\bibitem{kauf}
L.H. Kauffman, State Models and the Jones Polynomial, Topology, 26, 395-407, 1987. 

\bibitem{kim}
Peter Kim, Lee Stemkoski, Cornelia Yuen, Polynomial knots of
degree Five, MIT undergraduate Journal of Mathematics, 2001.

\bibitem{kui}
 Nicolas Kuiper, A new knot invariant, Mathematische, Annalen, 278, 193-209.
\bibitem{rm1}
R. Mishra, Polynomial representation of torus knots of type
{$(p,q)$}, J. Knot Theory Ramifications 8 (1999) 667-700.
\bibitem{rm2}
R. Mishra, Minimal degree sequence for torus knots, J. Knot Theory
Ramifications 9 (2000) 759-769.
\bibitem{rm3}
R.Mishra, Polynomial representation of strongly invertible and
strongly negative amphicheiral knots,  Osaka Journal of
Mathematics,Vol. 43, no.3, 2006.
\bibitem{p1}
P. Madeti and R. Mishra, Minimal Degree Sequence for Torus Knots of
Type $(p,2p-1)$, J. Knot Theory Ramifications,Vol. 15, No.9 (2006).

\bibitem{p2}
P. Madeti and R. Mishra, Minimal Degree Sequence for 2-bridge knots,
Fund. Math. 190 (2006) 191-210.
\bibitem{p3}
P. Madeti and R. Mishra, Minimal Degree Sequence for Torus Knots
of Type $(p, q)$, J. Knot Theory Ramifications, Vol. 18, No. 4
(2009), 485-491.

\bibitem{mcf}
McFeron, Donovan, The minimal degree sequence of the polynomial
figure eight knot (REU  2002)

\bibitem{dp2}
D. Pecker, Sur le th\'{e}or\`{e}me local de Harnack, C. R. Acad.
Sci. Paris, t. 326, S\'{e}rie 1, 1998, pp. 573-576.
\bibitem{rand}
R.Randel, Invariants of Piecewise Linear Knots, Knot theory, Banach Center Publications, Vol 42,Institute of Mathematics, Polish Academy of Sciences, Warszawa, 1998.

\bibitem{am}
R. Shukla and A. Ranjan, On Polynomial Representation of Torus
knots, J. of Knot Theory Ramifications 5 (1996) 279-294.
\bibitem{ars}
A. R. Shastri, Polynomial representations of knots, Tohoku Math.
J. (2) 44 (1992) 11-17.

\bibitem{schu}

 H. Schubert, Uber eine numerische Knoteninvariante,
 Math. Z. 61 (1954) 245–288.
 
\bibitem{rs}
R. Shukla, On polynomial isotopy of knot-types, Proc. Indian Acad.
Sci. Math. Sci. 104 (1994) 543-548.
\bibitem{trau}

A.K. Trautwein, An Introduction to Harmonic knots, Chapter 18, Ideal knots, Series on knots and everything, Vol 9, World Scientific (1998), 353-363. 



\bibitem{va}
Vassiliev, V.A., On Spaces of Polynomial knots, Sbornik:
Mathematics Volume 187 Number 2, 1996.

\end{thebibliography}
\end{document}